\documentclass[11pt,a4paper]{article}
\usepackage{url}
\usepackage[utf8]{inputenc}
\usepackage[english]{babel}
\usepackage{commath}
\usepackage[autostyle]{csquotes}
\usepackage{graphicx}
\graphicspath{{images/}}
\usepackage{enumerate}
\usepackage{float}
\usepackage{amsthm,amssymb,mathrsfs,setspace,amsmath}
\usepackage{caption}
\usepackage{subcaption}
\usepackage{algorithm}
\usepackage[shortlabels]{enumitem}
\usepackage[a4paper,bindingoffset=0.2in,left=0.5in,right=0.5in,top=1in,bottom=0.75in,footskip=.25in]{geometry}
\usepackage[colorlinks,citecolor=blue,urlcolor=blue,bookmarks=false,hypertexnames=true]{hyperref}
\numberwithin{figure}{subsection}
\newtheorem{theorem}{Theorem}[section]

\newtheorem{remark}{Remark}[section]

\newtheorem{note}{Note}
\newtheorem{lemma}{Lemma}
\DeclareUnicodeCharacter{2212}{-}

\newcommand{\A}{\alpha}
\newcommand{\Ph}{\phi}
\newcommand{\Lmd}{\lambda}
\newcommand{\LmdC}{\Lambda}
\newcommand{\M}{\mu}
\newcommand{\N}{\nu}

\newcommand{\G}{\gamma}
\newcommand{\GC}{\Gamma}
\newcommand{\X}{\xi}
\newcommand{\D}{\delta}
\newcommand{\DC}{\Delta}

\newcommand{\bigO}{\mathcal{O}}

\title{Generalized Non-Standard Finite Difference Method for Fractional PDEs on Non-Uniform Grids}
\author{Devank Mishra$^a$, Sheerin Kayenat$^b$, Amit K. Verma$^c$\thanks{$^c$akverma@iitp.ac.in (Corresponding Author)}\\\small{\it{$^{a}$ Department of Mathematics,}} \\\small{\it{Indian Institute of Technology Bhubaneshwar,}}\\\small{\it{Argul, Khordha, Odisha, India - 752050.}}\\\small{\it{$^{b}$Mathematics Division,}}\\ \small{\it{ School of Advanced Sciences \& Languages,}} \\\small{\it{VIT Bhopal University,}}\\\small{\it{ Bhopal, Madhya Pradesh, 466114, India.}}\\\small{\it{$^{c}$ Department of Mathematics,}} \\\small{\it{Indian Institute of Technology Patna,}}\\\small{\it{ Bihta, Patna 801106, (BR) India.}}}

\date{\today}

\begin{document}

\maketitle

\begin{abstract}
This paper proposes a novel Generalized Non-Standard Finite Difference (GNSFD) scheme for the numerical solution of a class of fractional partial differential equations (FrPDEs). The formulation of the method is grounded in optimization and leverages the fractional Taylor series (FrTS) expansion associated with Caputo fractional derivatives (FrDs). To discretize the  time derivatives, the non-trivial denominator functions are utilized. The theoretical analysis establishes the consistency, stability, and convergence of the proposed scheme. Results are compared against existing methods to substantiate the accuracy and computational efficiency of the approach.

\noindent {\it{Keywords:}} Least-squares; Fractional differential equations; Non-standard finite difference scheme; Caputo fractional derivative; Non-trivial denominator functions
\end{abstract}

\section{Introduction}
In recent years, Fractional Differential Equations (FrDEs) have garnered considerable attention in applied mathematics due to their ability to represent a wide range of physical processes. The significance and necessity of fractional calculus are demonstrated by numerous models and applications employing FrDs. They serve as a powerful framework for describing the memory effects and hereditary characteristics present in diverse materials and processes. Applications arise in areas such as signal processing, system identification, control engineering, and robotics \cite{baleanu2012fractional,zaslavsky2002chaos}, astrophysics and chemical physics \cite{li2007remarks,Miller1993fractional,Oldham1974fractional,Podlubny1998}, and other fields \cite{torvik1984appearance} are modeled using fractional differentiation and integration operators. Since exact solutions are often unavailable for FrDEs arising in physical applications, there is a need to design robust and precise computational techniques to effectively manage derivatives of arbitrary fractional order.  
Several types of FrDs are available in the literature \cite{TEODORO2019}. Among all, some popular non-local FrDs are the Riemann–Liouville derivative, the Caputo derivative, and the Grünwald–Letnikov derivative. \cite{Kilbas2006,Podlubny1998,Samko1987fractional,Usero2008}, etc. 

Over the last three decades, researchers have been trying to develop mesh-free methods (see \cite{benito2001influence, onate1996finite}). The generalized finite difference method (GFDM) is a strong form of mesh-free method, also referred to as mesh-free finite difference method (see \cite{GaveteMPE2015}). GFDM has been established as a reliable meshless technique for addressing several types of PDEs: advection–diffusion, wave propagation, beams, plates, etc. The method employed is a collocation-based meshless technique that integrates the Taylor expansion with the moving least square (LSQ) procedure designed to approximate the function as well as its derivative terms, while the point collocation strategy is used to construct the discrete algebraic system \cite{benito2001influence,gavete2003improvements,liszka1980finite}. In the GFDM, the moving LSQ procedure generates the function approximation using a distribution of nodes. The incorporation of expansion in a Taylor series enables a unified representation of the unknown function and its derivatives, thereby eliminating the need for the direct and computationally intensive evaluation of meshless shape function derivatives, as required in conventional moving LSQ-based approximations \cite{cheng2015meshless,li2016error,li2016stability,liu2009meshfree}. Through the point collocation scheme, the resulting algebraic system, sparse and banded in structure, is obtained by enforcing the governing partial DE at all interior nodes and applying the boundary conditions at the boundary nodes. Consequently, the GFDM avoids both mesh generation and numerical integration entirely \cite{benito2021solving,qu2021integrating,lin2021simulation,li2021meshless,benito2020solving,song2020generalized,gu2019generalized,
suchde2019meshfree}. 

Recently NSFD method has emerged as an alternative approach to resolve a diverse array of issues that involve algebraic, differential, biological, and chaotic systems in their mathematical models.
Comprehensive discussions on these methods can be found in the referenced books \cite{REM1994, REMb2000, REMb2005} and survey articles \cite{patidar2005use, patidar2016nonstandard}. The primary goal of these methods is to mitigate the numerical instabilities that arise in standard finite difference (SFD) approaches. In these methods, nontrivial denominator functions (DFs) such as $\sin h$ and $\left(e^{h} - 1\right)^2$ are used instead of standard denominator terms such as $h$ and $h^2$. Additionally, the non-linear term, for example, $(u_j^n)^2$, is represented by the product $u_{j-1}^{n+1} u_{j+1}^n$, which is non-local on the computational grid.

Hussian et al. \cite{hussian2008non} employed a non-standard discretization scheme to solve FrDEs. Moaddy et al. \cite{moaddy2011non} developed the NSFD scheme for addressing linear partial DEs involving FrDs in both space and time. Here, the FrDs were approximated using the Grünwald–Letnikov approach. Sweilam and Assiri \cite{sweilam2016non} developed NSFD with Riesz fractional definition to solve linear fractional hyperbolic partial DE. Zahra and Hikal \cite{zahra2017non} developed an algorithm based on the NSFD method to address a broad range of variable-order fractional optimal control problems formulated using either the Riemann–Liouville or Caputo derivatives. Wang et al. \cite{wang2020efficient} proposed an efficient NSFD scheme for numerically solving the chaotic fractional-order Chen system, where the discretization of the FrDE was carried out within a suitable nonlocal framework using the Grünwald–Letnikov approximation. Liu and Wang \cite{liu2021non} introduced an NSFD method for the space-fractional advection–diffusion equation and demonstrated, via the Fourier–Von Neumann stability analysis, that the scheme is unconditionally stable. Further, Taghipour and Aminikhah \cite{taghipour2022efficient} proposed an NSFD scheme that was well-posed, unconditionally stable, and convergent for the numerical solution of the distributed-order time-fractional reaction–diffusion equation. Here, the time-FrD is defined in the Caputo sense. 

In this paper, we present a novel non-standard version of GFDM. We refer to this method as the Generalised Non-standard Finite Difference Method. We have shown that if non-standard denominator functions are used concerning the time derivative, our results are very much improved. We also introduced another idea to use more than one denominator function in a given domain, which reduces the error significantly throughout the domain. To elaborate on the efficiency of the non-standard version of GFDM, we consider the following type of FrDE (\cite{Vargas2022FD}), based on Caputo's FrD (\cite{Caputo1967}):
\begin{eqnarray}\label{eq:mainProblem}
\begin{cases}
\displaystyle \mbox{FrDE}:~&u_t(t,x)= f_1(t,x,u)u_x(t,x)+f_2(t,x,u)u_{xx}(t,x)+g(t,x,u)D^\A_a(u(t,x))_x+ h(t,x,u),\\
\mbox{I.C.}:&u(0,x) = u_0(x),~x \in [0,1], \\
\mbox{B.C.}:&u(t,0) = u_L,\quad u(t,1) = u_R,~ t > 0,
\end{cases}
\end{eqnarray}
where  $t > 0$, $x \in [0,1]$, $D^\A_a(u(t,x))_x$ represents Caputo's FrD of $u$ w.r.t. $x$. The functions $f_1,f_2,g$ and $h$ are sufficiently smooth. The well-posedness of solutions in terms of existence, uniqueness, and regularity of \eqref{eq:mainProblem} has been discussed in \cite{baeumer2007fractional,Kilbas2006,Podlubny1998}. 

The paper is structured in seven sections. Section 2 deals with the finite difference approximation of the derivatives. In Section 3, we introduce a novel observation on non-trivial denominator functions and outline the algorithm for implementing the method. Section 4 establishes the assessment of the explicit scheme in terms of consistency, while Section 5 proves its convergence. Section 6 provides numerical experiments. Finally, Section 7 summarizes the main findings and conclusions.

\section{Discretization}
For the discretization of both integer and fractional-order spatial derivatives, we employ the GFDM, built upon the moving LSQ approach and utilizing the Taylor series expansion (see \cite{Vargas2022FD} and references therein).

The FrD term is discretized by using the Fractional Taylor Series (FrTS) expansion (see \cite{Usero2008}). Here, we use the following notation:
\begin{equation}
    D_a^{k\A}=D_a^{\A}\dots D_a^{\A}~(\text{$k$ times}).
\end{equation}
Certain regularity assumptions on the problem’s solution are required to ensure the validity of this expansion:
\begin{description}
\item[(A1)] For $k=1,2,3,~D_a^{k\A}u \in \mathcal{C}([0,1])$ and $D_a^{k\A}u \in \mathcal{I}_{\A} ([0,1])$,
\item[(A2)] $D_a^{4\A}u$ is continuous continuous over the interval $[0,1]$.
\end{description}
Following the notation in \cite{Usero2008}, $\mathcal{I}_{\alpha} ([0,1])$ is defined as
$$\mathcal{I}_{\alpha} ([0,1])=\{ f\in \mathcal{C}([0,1]):I_a^{\alpha}f(x)<\infty,~\forall x\in [0,1]\},$$ where $$I_a^{\alpha}f(x)=\frac{1}{\Gamma(\alpha)}\int_a^x (x-t)^{(\alpha-1)}f(t) dt.$$
We begin by discretizing the interval \([0,1]\) as 
$$M = \{0=x^{(0)}< x^{(1)}< \ldots< x^{(N+1)} = 1\},$$ consisting of \((N+2)\) nodes. For each node at which the discretization is to be performed, we select a subset \(S \subseteq M\) containing `$s$' nodes, referred to as a star. This subset is defined as $$ S = \{x_0, x_1, x_2, \ldots, x_s\} \subseteq M,$$ 
where \(x_0\) denotes the central node at which the discretization is computed. Various criteria for selecting the star corresponding to a central node are discussed in (\cite{Benito2008,Benito2003}).

\subsection{Discretization of spatial derivatives of classical (integer) order} \label{sec2.1}
For the approximation of derivatives of integral order appearing in equation \eqref{eq:mainProblem} at each node of the selected star, we employ a second-order Taylor series expansion, truncated after the quadratic term, of the $u(x_i)$ (as time $t$ is fixed, we are representing $u(t,x)$ simply by $u(x)$) about the central node $x_0$  
\begin{equation}
\label{eq:intTaylorExpansion}    u(x_i)=u(x_0)+\frac{\partial u(x_0)}{\partial x} (x_i-x_0)+\frac{1}{2}\frac{\partial^2 u(x_0)}{\partial x^2} (x_i-x_0)^2+\cdots.
\end{equation}
Let $U_i = u(x_i)$ denote the approximation to the continuous solution $U(x)$ at the node $x_i$. Furthermore, we define the vectors as 
\begin{equation*}
    d^T:=\begin{bmatrix}\frac{\partial U_i}{\partial x}& \frac{\partial^2 U_i}{\partial x^2}\end{bmatrix} = \begin{bmatrix}
        d_1 & d_2
    \end{bmatrix} ,\quad c_i^T:=\begin{bmatrix}h_i & \frac{h_i^2}{2} \end{bmatrix},
\end{equation*}
where $h_i=x_i-x_0$. Now, we define the operator $F$ as the sum of squares of weighted errors at each point in the star as:
\begin{equation*}
    \displaystyle F(d)=\sum_{i=1}^s (U_0-U_i+c_i^Td)^2w_i^2
\end{equation*}
where $w_i=w(x_0,x_i)$ is an appropriate weight function. Different properties and criteria for choosing weight functions can be found in \cite{lancaster1986curve}.
The minimization of $F$ concerning the vector $d$ is carried out as follows:
\begin{eqnarray}
&& \frac{\partial F}{\partial d}=0  \implies \begin{bmatrix}
       \displaystyle \frac{\partial F}{\partial d_1 } \\\displaystyle \frac{\partial F}{\partial d_2 }
    \end{bmatrix}=0,\\
&&    \begin{bmatrix}
      \displaystyle \sum_{i=1}^s 2(U_0-U_i+c_i^Td)w_i^2 h_i\\
      \displaystyle \sum_{i=1}^s 2(U_0-U_i+c_i^Td)w_i^2 \frac{h_i^2}{2}
    \end{bmatrix}=0,\\
&&  \displaystyle \sum_{i=1}^s w_i^2 \begin{bmatrix}
      h_i \\ 
      \displaystyle \frac{h_i^2}{2}
    \end{bmatrix}c_i^T d=\displaystyle -\sum_{i=1}^s (U_0-U_i)w_i^2\begin{bmatrix}
        h_i,\\ 
        \displaystyle \frac{h_i^2}{2}
    \end{bmatrix},\\
    &&\displaystyle \sum_{i=1}^s w_i^2 c_ic_i^T d=\displaystyle -\sum_{i=1}^s (U_0-U_i)w_i^2c_i.
\end{eqnarray}
As established in \cite{Gavete_2017}, the matrix $A=\displaystyle \sum_{i=1}^s w_i^2 c_ic_i^T$ is positive definite. Consequently, we obtain 
\begin{equation}\label{eq:intSystem}
    d=-U_0 \sum_{i=1}^s w_i^2 A^{-1}c_i + \sum_{i=1}^s U_iw_i^2 A^{-1}c_i.
\end{equation}
For convenience, we introduce the vectors
\begin{equation}
\label{eq:Lambda0} \Lmd_0:=\displaystyle \sum_{i=1}^s w_i^2 A^{-1}c_i=\begin{bmatrix}
    \Lmd_{0,1}\\ \Lmd_{0,2}
\end{bmatrix},~ 
\Lmd_i:=\displaystyle  w_i^2 A^{-1}c_i=\begin{bmatrix}
    \Lmd_{i,1}\\ \Lmd_{i,2}
\end{bmatrix}.
\end{equation}
for each $i \in \{1, \cdots, s\}$. From equations \eqref{eq:intSystem} and \eqref{eq:Lambda0}, we arrive at the final discretization of the integer order derivatives as below:
\begin{eqnarray}\label{eq:discretizationX}
    \begin{cases}
        \displaystyle \frac{\partial U_0}{\partial x}=-U_0\Lmd_{0,1} + \sum_{i=1}^s U_i\Lmd_{i,1}+ \bigO (h_i^2),\\
        \displaystyle \frac{\partial^2 U_0}{\partial x^2}=-U_0\Lmd_{0,2} + \sum_{i=1}^s U_i\Lmd_{i,2} + \bigO (h_i^2).
    \end{cases}
\end{eqnarray}

\subsection{Discretization of Caputo Derivative $\displaystyle D^{\A}_a u(t,x)$} \label{sec2.2}
The discretization of the Caputo derivative is obtained through the following procedure. Throughout, it is assumed that $u$ possesses sufficient smoothness to satisfy the conditions $(A1)$ and $(A2)$. Consider the FrTS expansion of a function $f(x)$ around $x_0$ as presented in \cite{Usero2008}:
    \begin{equation}
        f(x) = \sum_{j=0}^n D_a^{j\A}f(x_0)\frac{\DC_j(x)}{\GC (j\A+1)} + R_n(x,x_0,a),
    \end{equation}
    where  
    \begin{equation}
        \DC_k(x)=(x-a)^{k\A}-(x_0-a)^{k\A}+\GC (k\A+1) \sum_{j=1}^{k-1}\frac{(x_0-a)^{j\A} \DC_{k-j}(x)}{\GC (j\A+1) \GC ((k-j)\A +1)}.
    \end{equation}
Here, $R_n(x,x_0,a)$ denotes the remainder term of order $n \A$. 

Now, consider the FrTS of order $2 \A$ of a function $u$ at a nodal point $x_i$ centered at $x_0$ 
    \begin{equation}
        u(x_i) = u(x_0) + D_0^\A u(x_0) \frac{\DC_1(x_i)}{\GC(\A +1)} + D_0^{2\A} u(x_0) \frac{\DC_2(x_i)}{\GC(2\A +1)} + R_2(x,x_0,0),
    \end{equation}
    where
    \begin{equation}
        \DC_1(x_i)=x_i^\A - x_0^\A , \quad \DC_2(x_i) = x_i^{2\A} - x_0^{2\A} + \frac{\GC (2\A + 1)}{(\GC (\A +1 ))^2}x_o^\A \DC_1(x_i).
    \end{equation}
    We adopt the notation 
    \begin{align}
      &  \mathcal{D}^T=\begin{bmatrix}
            D_0^\A U_0 & D_0^{2\A}U_0
        \end{bmatrix},\quad C_i^T=\begin{bmatrix}
           \displaystyle \frac{\DC_1(x_i)}{\GC(\A +1)} & \displaystyle \frac{\DC_2(x_i)}{\GC(2\A +1)}
        \end{bmatrix},
    \end{align}
    We now define the weighted residual function as
    \begin{align}
       & F(\mathcal{D})=\sum_{i=1}^s (U_0-U_i+C_i^T \mathcal{D})^2w_i^2,
    \end{align}
 where $w_i=w(x_0,x_i)$ is an appropriate weight function. Minimizing $F$ with respect to $\mathcal{D}$ as in section \ref{sec2.1}, we arrive at the following system of linear equations:
    \begin{equation}
         \displaystyle \sum_{i=1}^s w_i^2 C_iC_i^T \mathcal{D}=\displaystyle -\sum_{i=1}^s (U_0-U_i)w_i^2C_i.
    \end{equation}
    We denote $\mathcal{A}:=\displaystyle \sum_{i=1}^s w_i^2 C_iC_i^T$, and arrive at the following lemma:
    \begin{lemma}
When conditions $(A1)$ and $(A2)$ are satisfied, the matrix $\mathcal{A}$ is positive definite.
   \end{lemma}
\renewcommand\qedsymbol{$\blacksquare$}
\begin{proof}
    For proof, see \cite{Vargas2022FD}. 
\end{proof}
\renewcommand\qedsymbol{QED}
\begin{remark}
If $\mathcal{A}$ is a positive definite matrix, then there exists a Cholesky decomposition $\mathcal{A}=LL^T$, where $L$ is a unique lower triangular matrix. This property ensures that ill-conditioned matrices remain absent during computation. For notational convenience, we continue to write $\mathcal{A}^{-1}$, while employing the Cholesky decomposition in the actual computations. 
\end{remark}
\begin{equation}\label{eq:caputoSystem}
    \mathcal{D}=-U_0 \sum_{i=1}^s w_i^2 \mathcal{A}^{-1}C_i + \sum_{i=1}^s U_iw_i^2 \mathcal{A}^{-1}C_i.
\end{equation}
For simplicity, if we denote \begin{equation}\label{eq:LambdaC0}\LmdC_0:=\displaystyle \sum_{i=1}^s w_i^2 \mathcal{A}^{-1}C_i=\begin{bmatrix}
    \LmdC_{0,1}\\ \LmdC_{0,2}
\end{bmatrix}~\mbox{and} ~ \LmdC_i:=\displaystyle  w_i^2 \mathcal{A}^{-1}C_i=\begin{bmatrix}
    \LmdC_{i,1}\\ \LmdC_{i,2}
\end{bmatrix},\end{equation}
for each $i \in \{1,2, \dots, s\}$. From equations \eqref{eq:caputoSystem} and \eqref{eq:LambdaC0}, we get our required discretization of fractional-order spatial derivatives:
\begin{eqnarray} \label{eq:discretizationC}
\begin{cases}
    \displaystyle D_0^\A U_0=-U_0\LmdC_{0,1} + \sum_{i=1}^s U_i\LmdC_{i,1}+ \bigO (\DC_2(x_i)), \\
    \displaystyle D_0^{2\A} U_0=-U_0\LmdC_{0,2} + \sum_{i=1}^s U_i\LmdC_{i,2}+ \bigO (\DC_2(x_i)).
\end{cases}
\end{eqnarray}

\subsection{Discretization of the temporal derivative of integer order}\label{sec2.3}
The time derivative is discretized using a first-order forward difference scheme as follows:
\begin{equation}\label{eq:discretizationT}
    \frac{\partial U(n\DC t,x_0)}{\partial t} = \frac{U((n+1)\DC t,x_0) - U(n\DC t,x_0)}{\DC t} + \bigO (\DC t).
\end{equation} 
\subsection{Numerical Scheme}
Substituting the discretizations from equations \eqref{eq:discretizationX}, \eqref{eq:discretizationC}, \eqref{eq:discretizationT} into equation \eqref{eq:mainProblem} yields the following numerical scheme: 
\begin{multline} \label{eq:standardNumericalScheme}
     \frac{U_0^{n+1}-U_0^n}{\DC t} = f_1(n\DC t, x_0 , U_0^n)(-U_0\Lmd_{0,1} + \sum_{i=1}^s U_i\Lmd_{i,1})+f_2(n\DC t, x_0 , U_0^n)(-U_0\Lmd_{0,2} + \sum_{i=1}^s U_i\Lmd_{i,2} ) \\ + g(n\DC t,x_0,U_0^n)(-U_0\LmdC_{0,1} + \sum_{i=1}^s U_i\LmdC_{i,1} ) +  h(n\DC t,x_0,U_0^n) + \bigO (\DC t , h_i^2 ,\DC_2(x_i)).
\end{multline}
\section{Generalized NSFD (GNSFD) Scheme}
In the non-standard technique, the term $\DC t$ in the denominator on the left-hand side of equation \eqref{eq:standardNumericalScheme} is replaced by a function $\Ph (\DC t)$ satisfying $\displaystyle \lim_{\DC t \to 0}\frac{\Ph (\DC t)}{\DC t}=1$. Consequently, we obtain
\begin{multline} \label{eq:nonStandardNumericalScheme}
     \frac{U_0^{n+1}-U_0^n}{\Ph (\DC (t))} = f_1(n\DC t, x_0 , U_0^n)(-U_0\Lmd_{0,1} + \sum_{i=1}^s U_i\Lmd_{i,1})+f_2(n\DC t, x_0 , U_0^n)(-U_0\Lmd_{0,2} + \sum_{i=1}^s U_i\Lmd_{i,2} ) \\ + g(n\DC t,x_0,U_0^n)(-U_0\LmdC_{0,1} + \sum_{i=1}^s U_i\LmdC_{i,1} ) +  h(n\DC t,x_0,U_0^n) + \bigO (\DC t , h_i^2 ,\DC_2(x_i))
\end{multline}
The selection of an appropriate denominator function, denoted by 
$\Ph (\DC t)$, can substantially improve the accuracy and stability of the solution. Several possible choices for 
$\Ph (\DC t)$ and their corresponding applications are discussed in \cite{kayenat2022choice}.

\subsection{Choice of Non-trivial Denominator Functions}
When a denominator function is employed for a given problem, the overall maximum error of the solution may not decrease uniformly. In practice, certain regions of the solution domain may exhibit substantial error reduction, while others show comparatively larger errors (see \cite{kayenat2022choice}). To address this, we propose to identify and use two denominator functions, $\Ph_1 (\DC t)$ and $\Ph_2 (\DC t)$, such that the regions of reduced error associated with each function are approximately complementary. We then construct a new denominator function as a linear combination, $$\A \Ph_1 (\DC t) + (1-\A) \Ph_2 (\DC t),$$ where $\A \in [0,1]$ is optimized locally at each node to minimize the error, thereby assigning greater weight to the function that results in lower errors at a particular node.

An additional strategy for improving accuracy is to introduce an extra parameter $\G$ into $\Ph (\DC t)$, resulting in $\Ph (\DC t,\G)$ with the condition $$\lim_{\DC t \to 0}\frac{\Ph (\DC (t),\G)}{\DC (t)}=1.$$ Optimizing $\G$ at each node further enhances the accuracy of the method. Combining these two ideas, we define 
\begin{equation*}
  \A \Ph_1 (\DC t,\G) + (1-\A) \Ph_2 (\DC t,\M)) 
\end{equation*} 
as the denominator function, where both $\A$ and the parameters $\G$ and $ \M$ are determined locally to minimize the error. This combined approach can yield significant improvements in the quality of the numerical solution.
\begin{note}
 This strategy of combining multiple denominator functions through a linear combination and introducing additional parameters can be extended beyond two functions. However, increasing the number of functions and parameters inevitably raises computational costs, and in some cases, the resulting improvement in accuracy may be marginal.
\end{note}
 \subsection{Algorithmic Implementation}
 We now outline an algorithm for implementing the proposed GNSFD method. Let the discretization of the interval $[0,1]$ be $\{0=x^{(0)}<x^{(1)}<\dots<x^{(N+1)}=1 \}$. The denominator function incorporates the parameters $\A$, $\G$ and $\M$, which are selected from the finite sets $\{\A_1, \A_2, \dots, \A_{n_1}\}$, $\{\G_1, \G_2, \dots, \G_{n_2}\}$, and $\{\M_1, \M_2, \dots, \M_{n_3}\}$ respectively.
Starting from the initial condition, the solution is computed at each node for the first time level \( t= n \Delta t \) (with \( n=0 \)), and subsequently for \( n=1\), This process is then repeated iteratively. At each time step, the solution at every node is available, and the values at the next time step \( t =(n+1) \Delta t \) are obtained according to the following procedure:
 \begin{algorithm}[H]
     \caption{Generalized Non-standard Finite Difference Scheme}
     \begin{enumerate}
         \item Set, minimumError $= +$ Infinity, $i=1$, $j=1$, $k=1$, $l=1$\\ $x_0=x^{(i)}, \A=\A_j, \G=\G_k, \M=\M_l$
         \item \label{calculation}Calculate $u(x_0,(n+1)\DC t)$ using equation \eqref{eq:nonStandardNumericalScheme}, and calculate currentError at the node.
         \item \textbf{IF:}currentError $<$ minimumError. \textbf{THEN:}Set, minimumError $=$ currentError and bestSolution $=u(x_0,(n+1)\DC t)$.
        \textbf{ELSE:GOTO:} Step \ref{changeParameter1}.
         \item \label{changeParameter1}\textbf{IF:}$j<n_1$ \textbf{THEN:}$j=j+1$, $\A=\A_j$ and
         \textbf{GOTO:} Step \ref{calculation}.
         \textbf{ELSE:}\textbf{GOTO:} Step \ref{changeParameter2}
         \item \label{changeParameter2} \textbf{IF:}$k<n_2$ \textbf{THEN:}$k=k+1$, $\G=\G_k$, $j=1,\A=\A_j$ and
         \textbf{GOTO:} Step \ref{calculation}.
         \textbf{ELSE:}\textbf{GOTO:} Step \ref{changeParameter3}.
         \item \label{changeParameter3} \textbf{IF:}$l<n_3$ \textbf{THEN:}$l=l+1$, $\M=\M_l$, $j=1,\A=\A_j$, $k=1,\G=\G_k$ and
         \textbf{GOTO:} Step \ref{calculation}.
         \textbf{ELSE:}\textbf{GOTO:} Step \ref{changeNode}.
         \item \label{changeNode} Set, final solution = bestSolution.
         \textbf{IF:}$i<N$, \textbf{THEN:}$i=i+1,x_0=x^{(i)},j=1,k=1,l=1,\A=\A_j,\G=\G_k,\M=\M_l$ and \textbf{GOTO:}Step \ref{calculation}.
         \textbf{ELSE:GOTO:} Step \ref{end}.
         \item \label{end} Required solution for time $(n+1)\DC t$ is obtained.
     \end{enumerate}
 \end{algorithm}
 \begin{note}
     The algorithm provided above for choosing best values for parameters $\A,\G,\M$ for the purpose of getting optimized value at a node may not be very efficient and a more efficient method can be developed for the same. 
 \end{note}  
\section{Consistency of the method}
The consistency of the proposed GNSFD method can be concluded by using the analysis for integer-order derivatives in \cite{Gavete_2017} and for non-integer-order Caputo derivatives in \cite{Vargas2022FD}. 

\section{Convergence of the Method}
The convergence of the proposed method follows directly from Theorem 4.1 in \cite{Vargas2022FD}, with minor modifications. We assume that all the hypotheses of Theorem 4.1 in \cite{Vargas2022FD} hold. Let the denominator function be defined as $\Ph(\DC t) = \A \Ph_1(\DC t, \G) + (1 - \A) \Ph_2(\DC t, \M)$, and let $u(n \DC t, x_0)$ denote the numerical approximation of $U_0^n$ satisfying equation \eqref{eq:nonStandardNumericalScheme}. Consequently, the exact solution $u_0^n$ will satisfy
\begin{multline} \label{eq:nonStandardNumericalSchemeReplacedWithSmallU}
     \frac{u_0^{n+1}-u_0^n}{\Ph (\DC (t))} = f_1(n\DC t, x_0 , u_0^n)(-u_0\Lmd_{0,1} + \sum_{i=1}^s u_i\Lmd_{i,1})+f_2(n\DC t, x_0 , u_0^n)(-u_0\Lmd_{0,2} + \sum_{i=1}^s u_i\Lmd_{i,2} ) \\ + g(n\DC t,x_0,u_0^n)(-u_0\LmdC_{0,1} + \sum_{i=1}^s u_i\LmdC_{i,1} ) +  h(n\DC t,x_0,u_0^n) + \bigO (\DC t , h_i^2 ,\DC_2(x_i)).
\end{multline} 
Subtracting equation \eqref{eq:nonStandardNumericalScheme} from equation \eqref{eq:nonStandardNumericalSchemeReplacedWithSmallU} and representing $e_i^n=u_i^n-U_i^n$, we get
\begin{multline}\label{eq:errorFirst}
     \frac{e_0^{n+1}-e_0^n}{\Ph (\DC (t))}=-\Lmd_{0,1}(f_{1,0}^nu_0^n-F_{1,0}^nU_0^n)+(f_{1,0}^n\sum_{i=1}^s\Lmd_{i,1}u_i^n-F_{1,0}^n\sum_{i=1}^s\Lmd_{i,1}U_i^n) \\ 
    -\Lmd_{0,2}(f_{2,0}^nu_0^n-F_{2,0}^nU_0^n)+(f_{2,0}^n\sum_{i=1}^s\Lmd_{i,2}u_i^n-F_{2,0}^n\sum_{i=1}^s\Lmd_{i,2}U_i^n)\\
    -\LmdC_{0,1}(g_1^nu_0^n-G_1^nU_0^n)+(g_1^n\sum_{i=1}^s\LmdC_{i,1}u_i^n-G_1^n\sum_{i=1}^s\LmdC_{i,1}U_i^n)\\+(h_0^n-H_0^n)+\bigO (\DC t , h_i^2 ,\DC_2(x_i)),
\end{multline}
where $$f_{i,0}^n=f_i(n\DC t,x_0,u_0^n), F_{i,0}^n=f_i(n\DC t,x_0,U_0^n)~ \text{for}~ i=1,2 ~\text{and} ~$$ 
$$ g_0^n=g(n\DC t,x_0,u_0^n), G_0^n=g(n\DC t,x_0,U_0^n), h_0^n=h(n\DC t,x_0,u_0^n), H_0^n=h(n\DC t,x_0,U_0^n).$$ 
Now, consider
\begin{align}
    f_{1,0}^nu_0^n-F_{1,0}^nU_0^n &=f_{1,0}^nu_0^n-f_{1,0}^nU_0^n+f_{1,0}^nU_0^n-F_{1,0}^nU_0^n\\
    &=f_{1,0}^n(u_0^n-U_0^n)+(f_{1,0}^n-F_{1,0}^n)U_0^n\\
    &=f_{1,0}^ne_0^n+(f_{1,0}^n-F_{1,0}^n)U_0^n\\
    \label{eq:X1}
    &=[f_{1,0}^n+\frac{(f_{1,0}^n-F_{1,0}^n)}{(u_0^n-U_0^n)}U_0^n]e_0^n.
\end{align}
Now using the Mean Value Theorem for some $\X_1 \in (u_0^n,U_0^n)\bigcup (U_0^n,u_0^n)$, we have
\begin{equation*}
    \frac{(f_{1,0}^n-F_{1,0}^n)}{(u_0^n-U_0^n)}=\frac{\partial f_1}{\partial u}(\X_1) ,\; \text{put it in equation\eqref{eq:X1} to get}
\end{equation*}
\begin{equation}
    f_{1,0}^nu_0^n-F_{1,0}^nU_0^n=[f_{1,0}^n+\frac{\partial f_1}{\partial u}(\X_1)U_0^n]e_0^n.
\end{equation}
Again, consider 
\begin{align*}    f_{1,0}^n\sum_{i=1}^s\Lmd_{i,1}u_i^n-F_{1,0}^n\sum_{i=1}^s\Lmd_{i,1}U_i^n &=f_{1,0}^n\sum_{i=1}^s\Lmd_{i,1}u_i^n-f_{1,0}^n\sum_{i=1}^s\Lmd_{i,1}U_i^n+f_{1,0}^n\sum_{i=1}^s\Lmd_{i,1}U_i^n-F_{1,0}^n\sum_{i=1}^s\Lmd_{i,1}U_i^n\\
&=f_{1,0}^n\sum_{i=1}^s\Lmd_{i,1}(u_i^n-U_i^n)+(f_{1,0}^n-F_{1,0}^n)\sum_{i=1}^s\Lmd_{i,1}U_i^n\\
&=f_{1,0}^n\sum_{i=1}^s\Lmd_{i,1}e_i^n+\frac{\partial f_1}{\partial u}(\X_1)e_0^n\sum_{i=1}^s\Lmd_{i,1}U_i^n.
\end{align*}
Similarly, we proceed for the terms having $f_2,g,h$ to get $\X_2,\X_3,\X_4$, respectively. Then, using all these in equation \eqref{eq:errorFirst}, we deduce
\begin{multline}
     \frac{e_0^{n+1}-e_0^n}{\Ph (\DC (t))}=-\Lmd_{0,1}[f_{1,0}^n+\frac{\partial f_1}{\partial u}(\X_1)U_0^n]e_0^n+f_{1,0}^n\sum_{i=1}^s\Lmd_{i,1}e_i^n+\frac{\partial f_1}{\partial u}(\X_1)e_0^n\sum_{i=1}^s\Lmd_{i,1}U_i^n\\
    -\Lmd_{0,2}[f_{2,0}^n+\frac{\partial f_2}{\partial u}(\X_2)U_0^n]e_0^n+f_{2,0}^n\sum_{i=1}^s\Lmd_{i,2}e_i^n+\frac{\partial f_2}{\partial u}(\X_2)e_0^n\sum_{i=1}^s\Lmd_{i,2}U_i^n\\
    -\LmdC_{0,1}[g_0^n+\frac{\partial g}{\partial u}(\X_3)U_0^n]e_0^n+g_0^n\sum_{i=1}^s\LmdC_{i,1}e_i^n+\frac{\partial g}{\partial u}(\X_3)e_0^n\sum_{i=1}^s\LmdC_{i,1}U_i^n\\
    \frac{\partial h}{\partial u}(\X_4)e_0^n+\bigO (\DC t, h_i^2,\DC_2(x_i)),\\
    e_0^{n+1}=e_0^n[1-\Ph(\DC t)(\Lmd_{0,1}[f_{1,0}^n+\frac{\partial f_1}{\partial u}(\X_1)U_0^n]+\Lmd_{0,2}[f_{2,0}^n+\frac{\partial f_2}{\partial u}(\X_2)U_0^n]+\LmdC_{0,1}[g_0^n+\frac{\partial g}{\partial u}(\X_3)U_0^n]\\-\frac{\partial f_1}{\partial u}(\X_1)\sum_{i=1}^s\Lmd_{i,1}U_i^n-\frac{\partial f_2}{\partial u}(\X_2)\sum_{i=1}^s\Lmd_{i,2}U_i^n-\frac{\partial g}{\partial u}(\X_3)\sum_{i=1}^s\LmdC_{i,1}U_i^n-\frac{\partial h}{\partial u}(\X_4))]\\
    +\Ph(\DC t)[f_{1,0}^n\sum_{i=1}^s\Lmd_{i,1}e_i^n+f_{2,0}^n\sum_{i=1}^s\Lmd_{i,2}e_i^n+g_0^n\sum_{i=1}^s\LmdC_{i,1}e_i^n]+\bigO (\DC t , h_i^2 ,\DC_2(x_i)).
\end{multline}
Formally, let $e^n:=\max_{i\in \{ 0,1,\dots,s \}}\{ |e_i^n| \}$. Utilizing this notation and applying it to bound the preceding expression, we obtain
\begin{equation}\label{eq:errorSecond}
    e^{n+1} \leq e^n \left[ \left| 1-\Ph(\DC t)\left( A\Lmd_{0,1}+B\Lmd_{0,2}+C\LmdC_{0,1}+D \right) \right| +\Ph(\DC t) E \right] +\bigO (\DC t , h_i^2 ,\DC_2(x_i)),
\end{equation}
where,
\begin{align*}
    &A=[f_{1,0}^n+\frac{\partial f_1}{\partial u}(\X_1)U_0^n],~ B=[f_{2,0}^n+\frac{\partial f_2}{\partial u}(\X_2)U_0^n], ~C=[g_0^n+\frac{\partial g}{\partial u}(\X_3)U_0^n], \\
    &D=\left( -\frac{\partial f_1}{\partial u}(\X_1)\sum_{i=1}^s\Lmd_{i,1}U_i^n-\frac{\partial f_2}{\partial u}(\X_2)\sum_{i=1}^s\Lmd_{i,2}U_i^n-\frac{\partial g}{\partial u}(\X_3)\sum_{i=1}^s\LmdC_{i,1}U_i^n-\frac{\partial h}{\partial u}(\X_4) \right), \\
    &E=\left( |f_{1,0}^n|\sum_{i=1}^s|\Lmd_{i,1}|+|f_{2,0}^n|\sum_{i=1}^s|\Lmd_{i,2}|+|g_0^n|\sum_{i=1}^s|\LmdC_{i,1}| \right).
\end{align*}
Hence from equation \eqref{eq:errorSecond}, it follows the condition for convergence which is given by
\begin{eqnarray}
    &&\left| 1-\Ph(\DC t)\left( A\Lmd_{0,1}+B\Lmd_{0,2}+C\LmdC_{0,1}+D \right) \right| +\Ph(\DC t) E < 1\\
    && -1+\Ph(\DC t) E < 1-\Ph(\DC t)\left( A\Lmd_{0,1}+B\Lmd_{0,2}+C\LmdC_{0,1}+D \right) < 1-\Ph(\DC t) E\\
  && \text{First inequality:}  -1+\Ph(\DC t) E < 1-\Ph(\DC t)\left( A\Lmd_{0,1}+B\Lmd_{0,2}+C\LmdC_{0,1}+D \right)\\
&& \label{eq:condition2}
    \Longrightarrow \Ph(\DC t) < \frac{2}{A\Lmd_{0,1}+B\Lmd_{0,2}+C\LmdC_{0,1}+D+E}.\\
    && \text{Second inequality:}  1-\Ph(\DC t)\left( A\Lmd_{0,1}+B\Lmd_{0,2}+C\LmdC_{0,1}+D \right) < 1-\Ph(\DC t) E\\
&& \label{eq:condition1}
    \Longrightarrow A\Lmd_{0,1}+B\Lmd_{0,2}+C\LmdC_{0,1}+D-E > 0
\end{eqnarray}
As our $\Ph(\DC t),\Lmd_{0,1},\Lmd_{0,2},\LmdC_{0,1}$ depends on parameters $\A,\G,\M$ also, so our conditions  \eqref{eq:condition2} , \eqref{eq:condition1} should be satisfied for at least one value of these parameters.
\begin{theorem}\cite{Vargas2022FD}
Let $u(t,x)$ satisfies (A1) and (A2) and $A\Lmd_{0,1}+B\Lmd_{0,2}+C\LmdC_{0,1}+D+E>0$. Further let that $f_1$, $f_2$, $g$ and $h$ be differentiable w.r.t $u$. 
Under the previous notations, let $U_0^n$ represent $u$ at $(n\Delta t, x_0)$. Als
Then, the explicit scheme \eqref{eq:nonStandardNumericalScheme} is convergent if
\[
\Ph(\DC t) < \frac{2}{A\Lmd_{0,1}+B\Lmd_{0,2}+C\LmdC_{0,1}+D+E}\;  \text{and}\;\; A\Lmd_{0,1}+B\Lmd_{0,2}+C\LmdC_{0,1}+D-E > 0\] 
for some constants $A, B, C, D$ and $E$ depending on the distribution of the nodes and explicitly computed in the proof.
\end{theorem}

\section{Numerical Illustrations}
This section presents examples, each satisfying the conditions stated in Theorem 4.1, to illustrate the performance of the proposed method in terms of both efficiency and accuracy. The computations are carried out on irregular node distributions. 

\subsection{Example 1} 
We consider the following FrPDE  
 \begin{equation}
 \frac{\partial u(t,x)}{\partial t}=\GC (1.2)x^{1.8}D^{1.8}_0 u(t,x)+(6x^3-3x^2)e^{-t},\quad x \in [0,1],\; t>0,
 \end{equation} 
with initial condition $u(0,x)=x^2-x^3$ and boundary conditions $u(t,0)=0$, $u(t,1)=0$. 
It's exact solution is given by $$u(t,x)=(x^2-x^3)e^{-t}.$$

\begin{figure}[H]
    \centering
    \includegraphics[height=7cm, width=9cm]{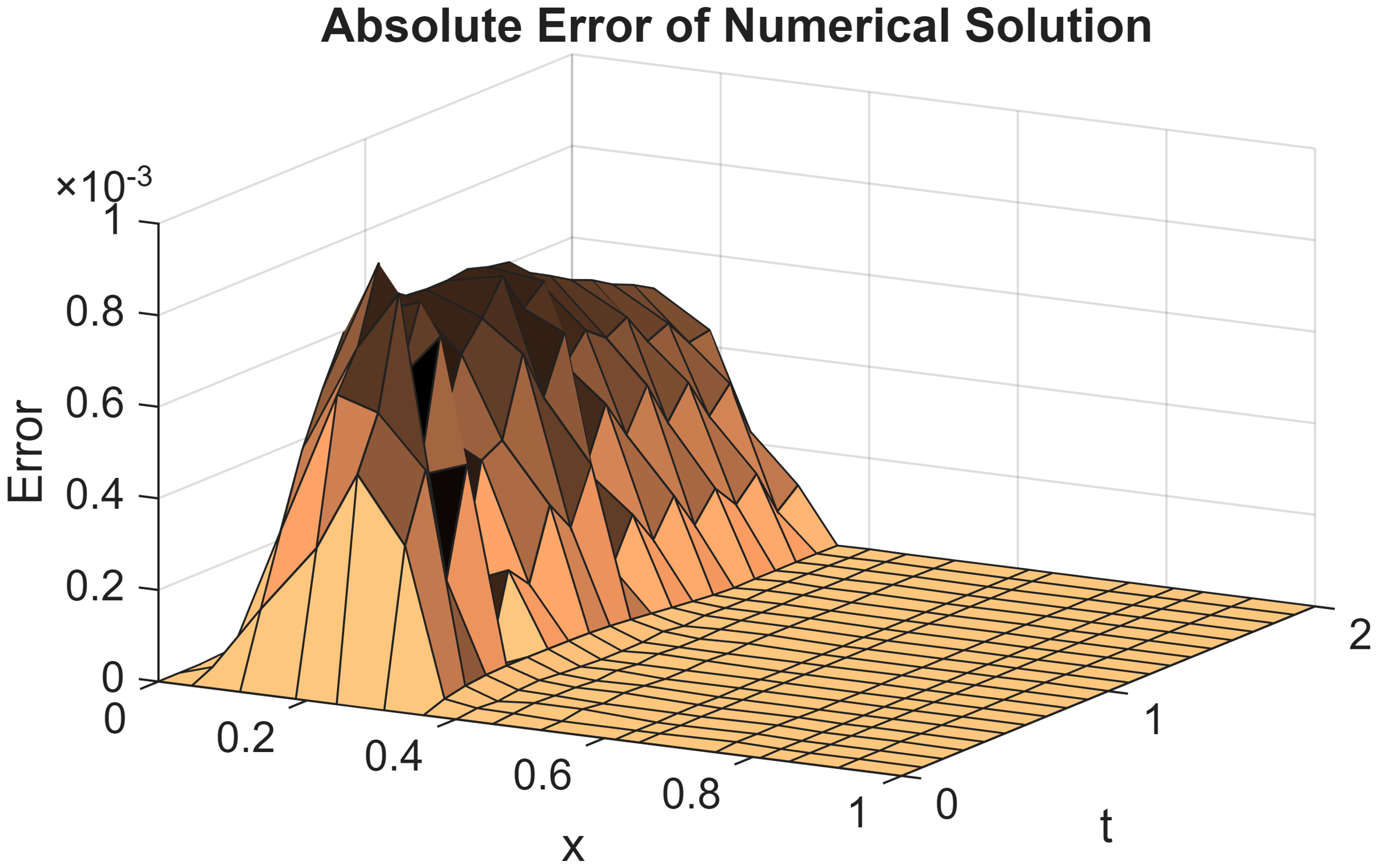}
    \caption{Absolute error at each grid point for $\Ph_1(\DC t,\G) $ as denominator function, for example 1 }
    \label{fig:3dErrorEx1Exponential}
\end{figure}

\begin{figure}[H]
    \centering
    \includegraphics[height=7cm, width=9cm]{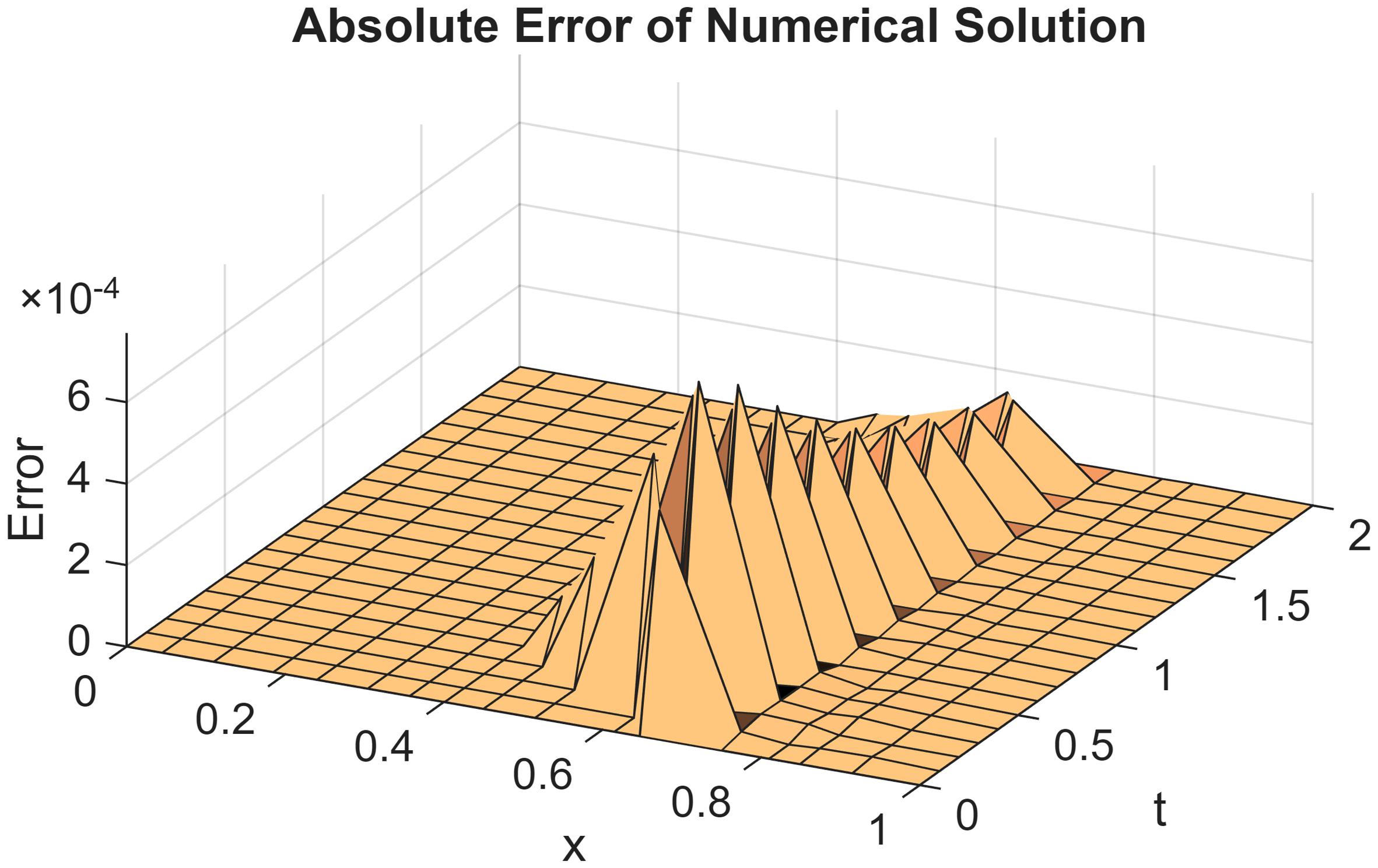}
    \caption{Absolute error at each grid point for example 1, taking $\Ph_2(\DC t,\M) $ as denominator function}
    \label{fig:3dErrorEx1Sin}
\end{figure}

 We have taken irregular discretization of $[0,1]$ with 17 nodes and $\DC t = 0.1$ as 
 $$\{ 0,0.045,0.11,0.185,0.24,0.304,0.357,0.401,0.45,0.515,0.615,0.647,0.75,0.81,0.849,0.915,1\}.$$
Taking $\displaystyle \Ph_1(\DC t,\G)=\frac{(e^{\G k}-1)}{\G}$ as denominator function of the non-standard scheme described in the equation \eqref{eq:nonStandardNumericalScheme} we get relatively higher errors for some initial nodes as compared to remaining nodes. This becomes clearer from the figure \ref{fig:3dErrorEx1Exponential}.
\paragraph{}
Taking $\displaystyle \Ph_2(\DC t,\M)=\frac{\sin(\M k)}{\M}$ as a denominator function, we get relatively higher errors for some middle nodes as compared to remaining (see \ref{fig:3dErrorEx1Sin}).
\paragraph{}
These observations are enough to motivate us towards the novel idea of using a denominator function that is a linear combination of these two. At each node, we will vary the weights so that the function that produces less error will be given more weight.
\paragraph{}

Now we take $\displaystyle \Ph_1(\DC t,\G)=\frac{(e^{\G k}-1)}{\G},\;\Ph_2(\DC t,\M)=\frac{\sin(\M k)}{\M}$, where $\G,\M\in\{1,2,\dots,10 \}$ and arrive at our final denominator function given as $$\A \Ph_1 (\DC t,\G) + (1-\A) \Ph_2 (\DC t,\M)$$ where $\A\in\{0,0.1,\dots,1.0 \}$. 

\begin{figure}[H]
    \centering
\includegraphics[height=7cm, width=9cm]{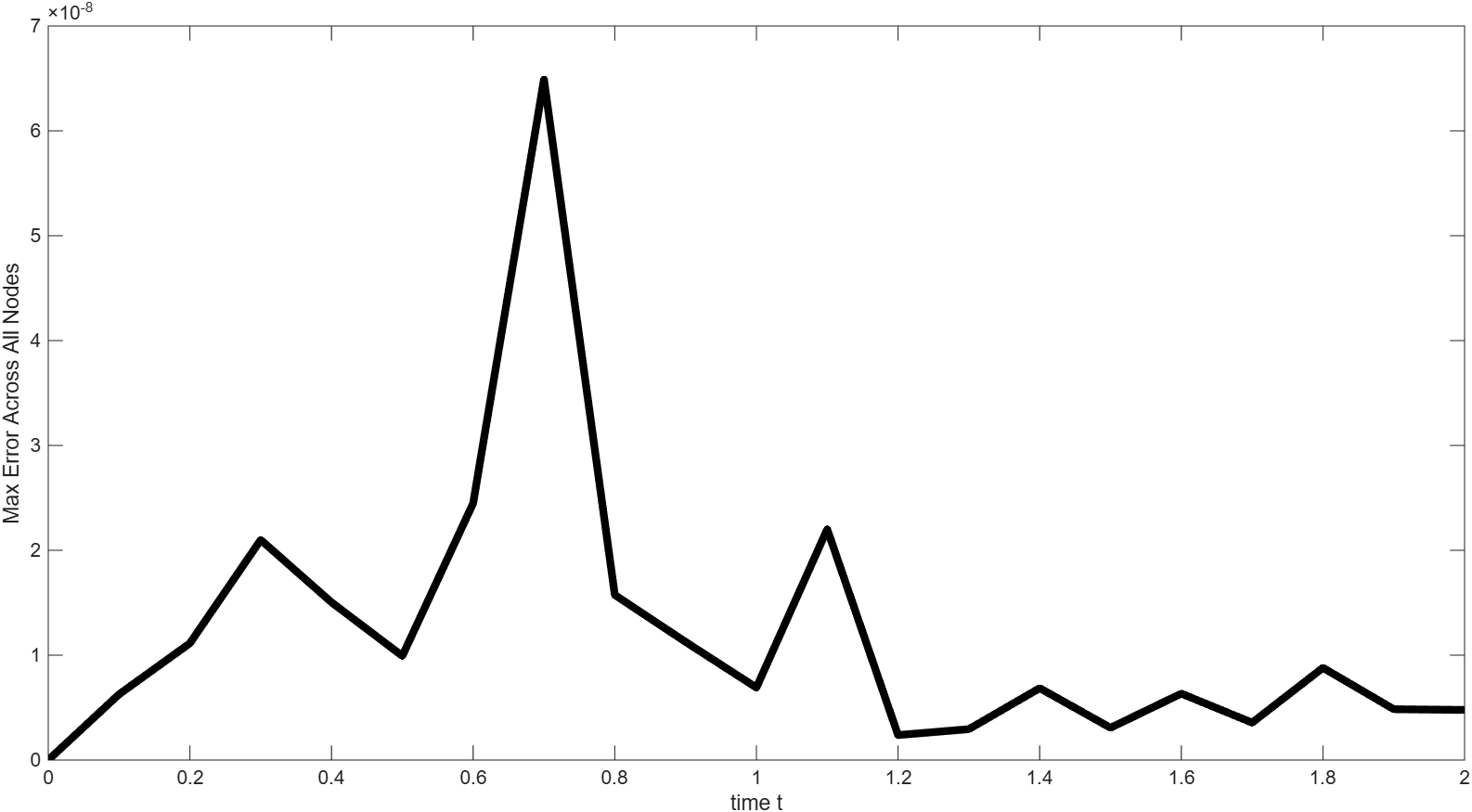}
    \caption{Max error at different values of t for example 1 on taking $(\A \Ph_1 (\DC t,\G) + (1-\A) \Ph_2 (\DC t,\M))$ as denominator function}
    \label{fig:2dErrorEx1}
\end{figure}

\begin{figure}[H]
    \centering
    \includegraphics[height=7cm, width=9cm]{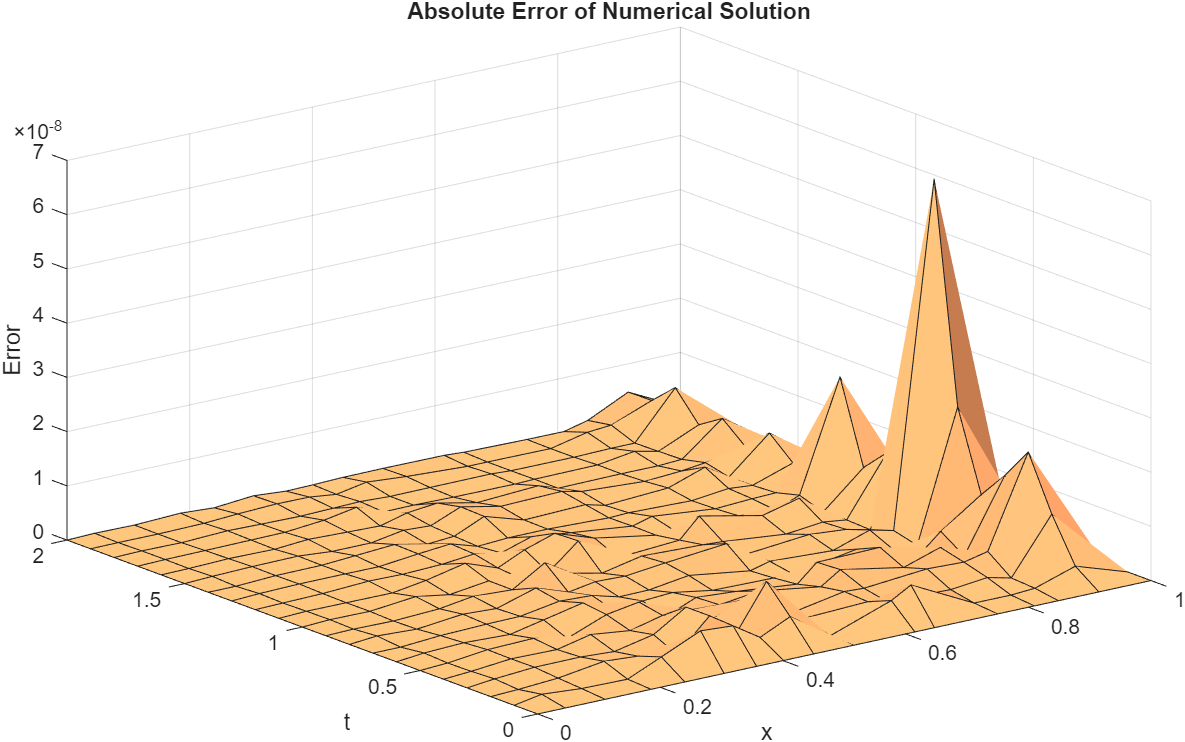}
    \caption{Absolute error at each grid point for example 1, on taking $(\A \Ph_1 (\DC t,\G) + (1-\A) \Ph_2 (\DC t,\M))$ as denominator function}
    \label{fig:3dErrorEx1}
\end{figure}

For $\DC t=0.1$, we compare the errors produced by this method with those given in
\begin{center}
    \begin{tabular}{||c|c|c||}
    \hline
       Time Step  & Max Error & Errors in \cite{Vargas2022FD} \\
       \hline \hline
        0.5 & $9.9276\times 10^{-9}$ & $8.3543\times 10^{-3}$ \\
        \hline
        1.0 & $6.9068\times 10^{-9}$ & $6.7542\times 10^{-3}$ \\
        \hline
        1.5 & $3.0760\times 10^{-9}$ & $9.9786\times 10^{-4}$ \\
        \hline
        2.0 & $4.7791\times 10^{-9}$ & $3.5349\times 10^{-4}$ \\
        \hline
    \end{tabular}
\end{center}
Figures \ref{fig:2dErrorEx1}, \ref{fig:3dErrorEx1} depict the advantage of this approach, where error is suppressed throughout the domain except where the order is less than $10^{-8}$.
From the above figures and tables, it is verified that the use of more than one denominator function helps us reduce the error significantly throughout the domain. 
\subsection{Example 2} We consider the following PDEs
\begin{equation}
    \displaystyle \frac{\partial u}{\partial t} + u^{\M} \frac{\partial u}{\partial x} = \N u^\D \frac{\partial^2 u}{\partial x^2}
\end{equation}
where $\M = 1,\; \D = 1,\; \N=1$, with initial condition $ w(0,x)=1-e^{\frac{x}{\N}}+x(e^{\frac{1}{\N}}-1)$ (see \cite{VERMA2021APNUM}) and boundary conditions $$ u(t,0)=u(t,1)=0. $$
The exact solution is given by $ w(t,x)=\frac{1-e^{\frac{x}{\N}}+x(e^{\frac{1}{\N}}-1)}{1+t(e^{\frac{1}{\N}}-1)} $.

\begin{figure}[H]
    \centering
    \includegraphics[height=7cm, width=9cm]{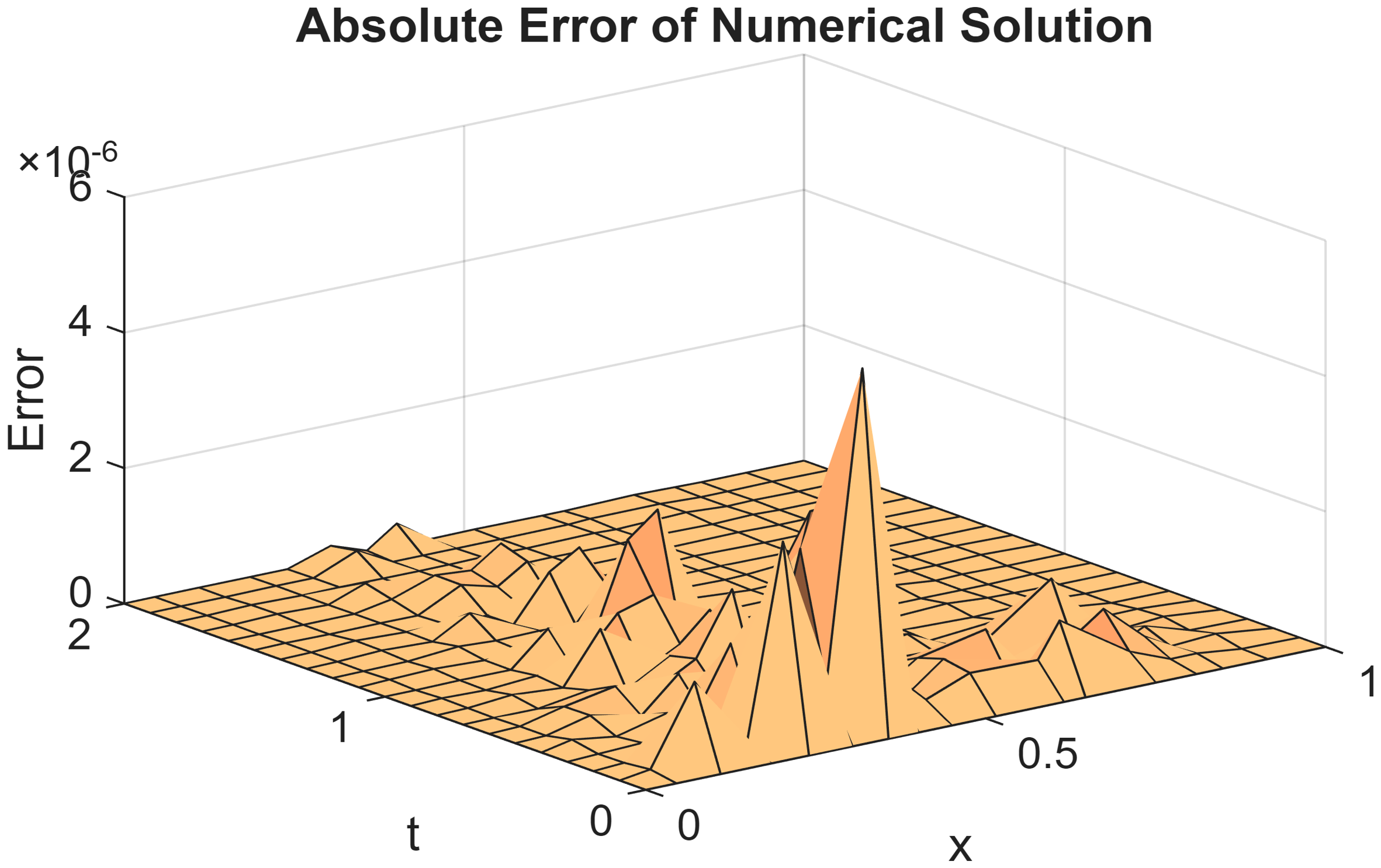}
    \caption{Absolute error at each grid point for example 2, taking $\Ph_1(\DC t,\G) $ as denominator function}
    \label{fig:3dErrorEx2Sin}
\end{figure}

\begin{figure}[H]
    \centering
    \includegraphics[height=7cm, width=9cm]{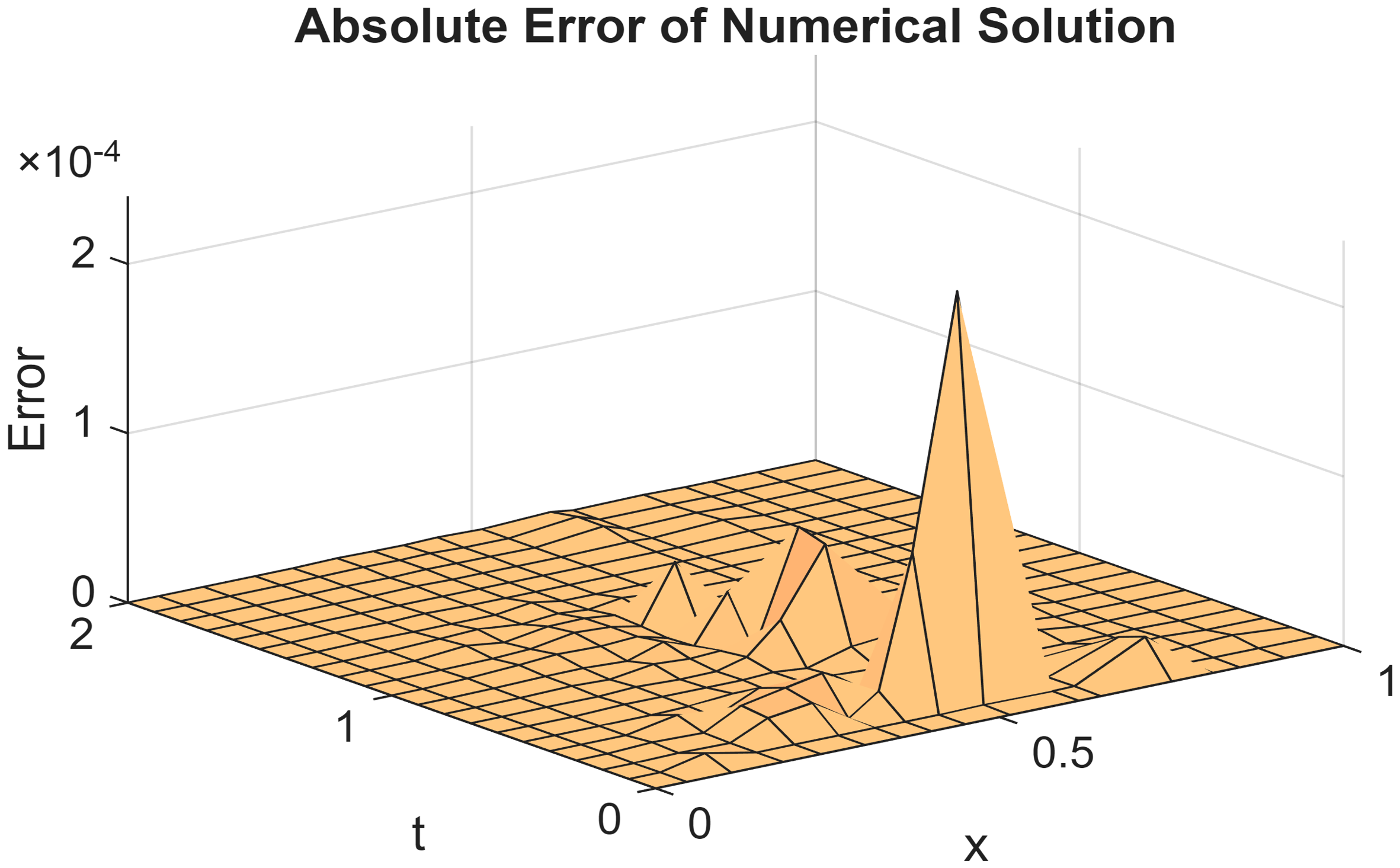}
    \caption{Absolute error at each grid point for example 2, taking $\Ph_2(\DC t,\theta) $ as denominator function}
    \label{fig:3dErrorEx2Tan}
\end{figure}

\begin{figure}[H]
    \centering
    \includegraphics[height=8cm, width=10cm]{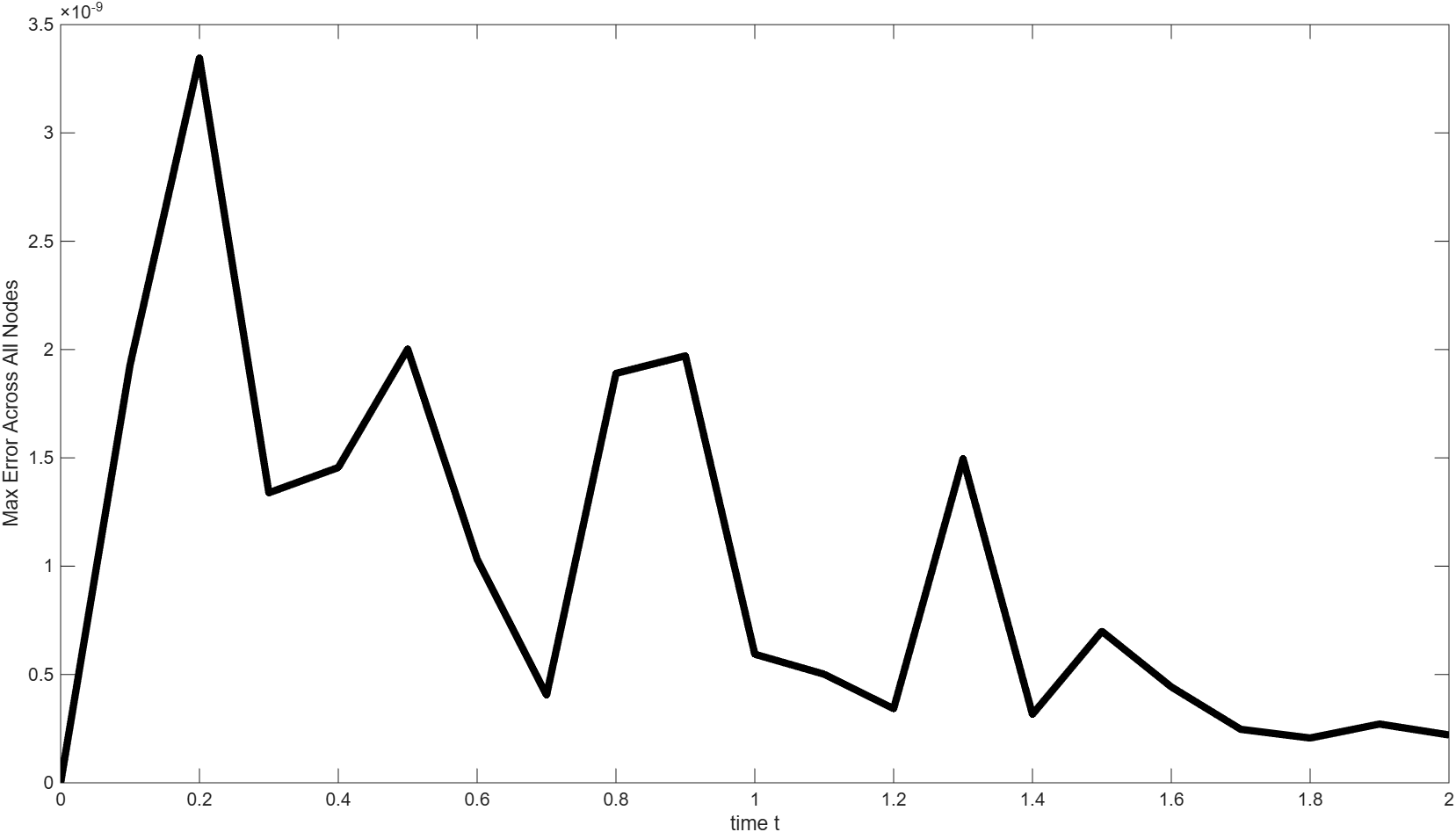}
    \caption{Max error at different values of $t$, for example 2 on taking $(\A \Ph_1 (\DC t,\G) + (1-\A) \Ph_2 (\DC t,\theta))$ as denominator function}
    \label{fig:2dErrorEx2}
\end{figure}

\begin{figure}[H]
    \centering
    \includegraphics[height=7cm, width=9cm]{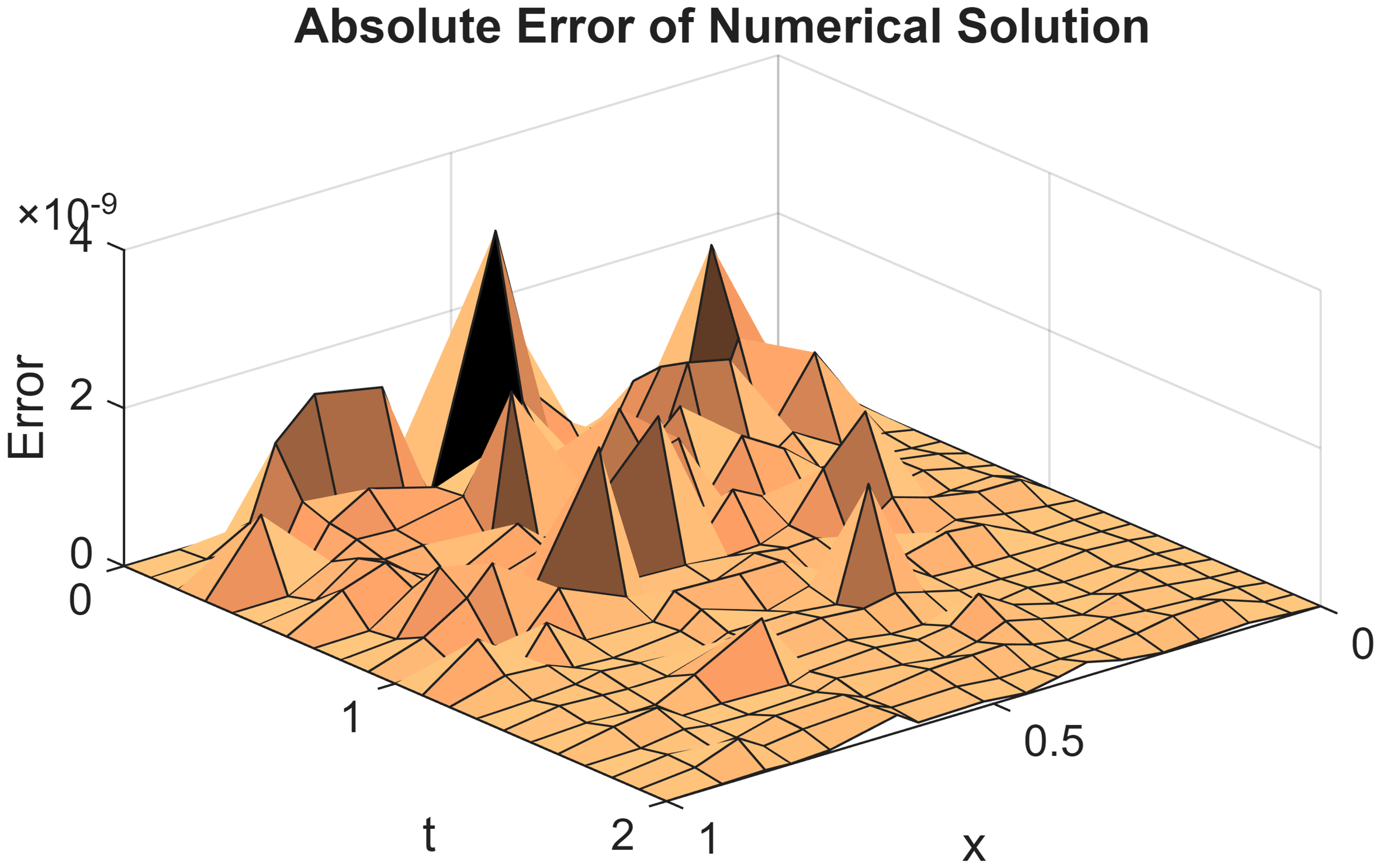}
    \caption{Absolute error at each grid point for example 2, on taking $(\A \Ph_1 (\DC t,\G) + (1-\A) \Ph_2 (\DC t,\theta))$ as denominator function}
    \label{fig:3dErrorEx2}
\end{figure}

\paragraph{}Using the same discretization as in the last example and $\DC t=0.1$ .
\paragraph{}
Taking $\displaystyle \Ph_1(\DC t,\G)=\frac{\sin(\G k)}{\G}$ as the denominator function (where $\G \in \{0.5,1,1.5,\dots,10.0\}$), Figure \ref{fig:3dErrorEx2Sin} presents error at each point on the grid. 
Taking $\displaystyle \Ph_2(\DC t,\theta)=\tan(\DC t)-\theta \sin^2 (\DC t)$ as the denominator function (where $\theta \in \{0,0.5,1,1.5,\dots,10.0\}$), Figure \ref{fig:3dErrorEx2Tan} presents error at each point on the grid. 
Now, similar to the last example, taking $\A \Ph_1 (\DC t,\G) + (1-\A) \Ph_2 (\DC t,\theta)$ as the denominator function where $\A\in\{0,0.1,\dots,1.0 \}$ we get the following, Figure \ref{fig:3dErrorEx2}, Figure \ref{fig:2dErrorEx2}
From the above figures and tables, it is verified that the use of more than one denominator function helps us reduce the error significantly throughout the domain. 
\section{Conclusion}
In this work, we have introduced a Generalized Non-Standard Finite Difference (GNSFD) scheme for the numerical solution of a class of fractional partial differential equations (FrPDEs). The scheme is constructed using the fractional Taylor series expansion in conjunction with Caputo fractional derivatives, and incorporates non-trivial denominator functions to discretize the fractional time derivatives. We have shown that instead of focusing on one denominator function for the entire domain, if we try a combination of denominator functions, it may provide even more accurate results. This fact has been validated by two test examples.   

\section*{Conflict of interest}
The authors declare no conflict of interest. 

\section*{Data availability} Data sharing does not apply to this article, as no datasets were generated or analyzed during the current study.

\section*{Author contributions} All authors contributed equally to this work.

\section*{Acknowledgments} The first author is thankful to Dr. Amit Kumar Verma , Associate Professor, Department of Mathematics, IIT Patna, for providing all the necessary help and support for staying in IIT Patna during the summer internship. He is also thankful to his mother Mrs. Mani Mishra , his father Mr. Sunil Kumar Mishra and his friends Prabal , Anu for their constant support.

\bibliographystyle{plain}
\bibliography{projectAkv.bib}
\end{document}